\documentclass[11pt,reqno,a4]{amsart}
\usepackage{a4}
\usepackage{graphicx}
\usepackage{amsmath,amssymb,natbib,color,bbm,ifthen,graphicx,epsf}
\usepackage{hyperref}
\usepackage[latin1]{inputenc}

\newcommand{\eqsp}{\;}

\newtheorem{theorem}{Theorem}
\newtheorem{corollary}[theorem]{Corollary}
\newtheorem{lemma}[theorem]{Lemma}
\newtheorem{proposition}[theorem]{Proposition}

\theoremstyle{definition}

\newtheorem{hyp}{\text{A}}
\numberwithin{equation}{section}
\theoremstyle{remark}

\newcommand{\CPE}[3][]
{\ifthenelse{\equal{#1}{}}{\mathbbm{E}\left[\left. #2 \, \right| #3 \right]}{\operatorname{#1}\left[\left. #2 \, \right | #3 \right]}}
\newcommand{\CPP}[3][]
{\ifthenelse{\equal{#1}{}}{\mathbbm{P}\left(\left. #2 \, \right| #3 \right)}{\operatorname{#1}\left(\left. #2 \, \right | #3 \right)}}

\newcommand{\rhs}{right hand side}
\newcommand{\iid}{i.i.d.}
\newcommand{\wrt}{w.r.t.}

\newcommand{\PP}{\mathbb{P}}

\def\rme{\mathrm{e}}

\def\GG{\mathcal G}

\newcommand{\1}{\ensuremath{\mathbbm{1}}}

\def\<{\langle}
\def\>{\rangle}

\newcommand\pscal[2]{\langle #1, #2 \rangle}
\newcommand{\eqdef}{\ensuremath{\stackrel{\mathrm{def}}{=}}}

\def\Xset{\mathsf{X}}

\def\Xsigma{\mathcal{X}}

\def\ie{\textit{i.e.}}

\newcommand{\nset}{\mathbb{N}}
\newcommand{\rset}{\mathbb{R}}

\newcommand{\RK}{Q} 

\def\PP{\mathbb{P}}

\def\PE{\mathbb{E}}

\newcommand{\stdlimit}{\mbox{$\boldsymbol{\sigma}$}}

\def\cP{\check{P}}
\def\cPP{\check{\PP}}
\def\cPE{\check{\PE}}
\def\tPE{\tilde{\PE}}
\def\tPP{\tilde{\PP}}
\def\csigma{\check{\sigma}}
\def\ctau{\check{\tau}}
\def\cmu{\check{\mu}}

\newcommand\pfunc[1]{\mathcal{F}^+(#1)}
\newcommand{\Hphi}[1][]
{\ifthenelse{\equal{#1}{}}{\Phi}{\Phi_{#1}}}

\title[Limit Theorems for Subgeometric Markov Chains]{Bounds on Regeneration Times and Limit Theorems for Subgeometric Markov Chains\\
\vspace{6mm}
Bornes des temps de r\'eg\'en\'eration et th\'eor\`emes limites pour des cha\^ines de Markov sous-g\'eom\'etriques}

 \author{Randal Douc $^\star$}
 \address{Randal DOUC (Corresponding Author), CMAP, Ecole Polytechnique, 91128 Palaiseau Cedex, France}
 \email{douc@cmapx.polytechnique.fr}
 \author{Arnaud Guillin}
 \address{Arnaud GUILLIN, CEREMADE, Université Paris IX Dauphine, place du Mar\'echal
 de Lattre de Tassigny 75775 Paris cedex 16}
 \email{guillin@ceremade.dauphine.fr}
 \author{Eric Moulines}
 \address{Eric MOULINES, Ecole Nationale Supérieure des Télécommunications, Département de Traitement du Signal et des Images, 46, rue Barrault, 75634 PARIS Cédex 13, France}
 \email{moulines@tsi.enst.fr}

\begin{document}
\maketitle

\begin{abstract}
This paper studies limit theorems for Markov Chains with general state space under conditions which imply subgeometric ergodicity.
We obtain a central limit theorem and  moderate deviation principles  for additive not necessarily bounded functional of the Markov chains
under drift and minorization conditions which are weaker than the Foster-Lyapunov conditions.
The regeneration-split chain method and a precise control of the modulated moment of the hitting time to small sets
are employed in the proof.
\end{abstract}


\renewcommand{\abstractname}{Résumé}
\begin{abstract}
Nous \'etablissons dans ce papier des th\'eor\`emes limites pour des cha\^ines de Markov \`a espace d'\'etat g\'en\'eral sous des conditions impliquant l'ergodicit\'e sous g\'eom\'etrique. Sous des conditions de d\'erive et de minorisation plus faibles que celles de Foster-Lyapounov, nous obtenons un th\'eor\`eme de limite centrale et un principe de d\'eviation mod\'er\'ee pour des fonctionnelles additives non n\'ecessairement born\'ees de la cha\^{i}ne de Markov. La preuve repose sur la m\'ethode de r\'eg\'en\'eration et un contr\^ole pr\'ecis du moment modul\'e de temps d'atteinte d'ensembles petits.   
\end{abstract}
\newpage

%
%
%
%

\section{Introduction}
This paper studies  limit theorems and deviation inequalities  for a positive Harris recurrent Markov chain $\{X_k \}_{k \geq 0}$ on  a
general state space $\Xset$ equipped with a countably generated $\sigma$-field $\Xsigma$.
Results of this type for geometrically ergodic Markov chains are now well established:
see for instance \cite[Chapter 17]{meyn:tweedie:1993} for the central limit theorem and the law of iterated logarithm, \cite{acosta:chen:1998},
\cite{chen:1999}  for moderate deviations results.
However, the more subtle subgeometrical case is not nearly as well understood (see for instance \cite{djellout:guillin:2001}).

These results can be obtained by using the regeneration method constructed via the  splitting technique on returns to small sets. These
methods typically require bounds for modulated moments of the excursions between two regenerations. In practice, one most often control the
corresponding modulated moment of the excursion between two small set return times rather than regeneration times. Our first result in section
\ref{sec:boundsregenerationtime} relate these two bounds, extending to subgeometrical case results reported earlier in the geometric case
by \cite{roberts:tweedie:1999}. We then apply these bounds in sections \ref{sec:CLT-BerryEsseen}, \ref{sec:MDP} and \ref{sec:Deviation-Inequalities}.
In section \ref{sec:CLT-BerryEsseen}, we establish a CLT and Berry-Esseen bounds, sharpening estimates given in \cite{bolthausen:1982}. In section \ref{sec:MDP}, we establish a Moderate Deviation Principle for possibly
unbounded additive functionals of the Markov chains, extending results obtained earlier for bounded functionals and \emph{atomic} chains by \cite{djellout:guillin:2001}.
Finally, in section \ref{sec:Deviation-Inequalities}, we give deviation inequality for unbounded additive functionals of the Markov Chain.

Following \cite{nummelin:tuominen:1983}, we denote by $\Lambda_0$ the set of sequences such that
$r(n)$ is non decreasing and $\log r(n)/n \downarrow 0$ as $n \to \infty$ and by $\Lambda$ the set of
sequences for which $r(n) > 0$ for all $n \in \nset$ and for which there exists an $r_0 \in \Lambda_0$ which is equivalent to $r$
in the sense that
\[
0 < \liminf_{n \to \infty} \frac{r(n)}{r_0(n)} \quad \text{and} \quad \limsup_{n \to \infty} \frac{r(n)}{r_0(n)} < \infty \eqsp.
\]
Without loss of generality, we assume that $r(0)=1$ whenever $r \in \Lambda_0$.
Examples of subgeometric sequences include: polynomial sequences $r(n) = (n+1)^\delta$ ($\delta > 0$), or  subexponential sequences,
$r(n) = (n+1)^\delta e^{c n^\gamma}$ ($\delta > 0$, $c > 0$ and $\gamma \in (0,1)$).

Denote by $P$ the transition kernel of the chain and for $n \geq 1$, $P^n$ the $n$-th iterate of
the kernel.
For any signed measure $\mu$ on $(\Xset,\Xsigma)$, we denote by $\| \mu \|_f \eqdef \sup_{|g| \leq f} | \mu(g)|$ the $f$-total variation norm.
Let $f: \Xset \to [1,\infty)$ be a measurable function and $\{r(k)\} \in \Lambda$.
We shall call $\{ X_k \}$ $(f,r)$-ergodic (or $f$-ergodic at rate $\{r(k)\}$) if $P$ is aperiodic, $\phi$-irreducible and positive Harris recurrent Markov chain and
\begin{equation}
\label{eq:definition-(f,r)-ergodicity}
\lim_{n \to \infty} r(n)  \| P^n(x,\cdot) - \pi  \|_{f} < \infty \eqsp, \quad \text{for all}\ x \in \Xset \eqsp.
\end{equation}
where $\pi$ is the unique stationary distribution of the chain. If \eqref{eq:definition-(f,r)-ergodicity} holds for $f \equiv 1$,
then we call $\{X_k\}$ $r$-ergodic (or ergodic at rate $r$).
For positive Harris recurrent Markov chain \citep[Chapter V]{meyn:tweedie:1993} there exists some (and indeed infinitely many) small sets satisfying
for some constant $m$ and some probability measure $\nu$, the minorisation condition: $P^m(x,\cdot) \geq \epsilon \nu(\cdot)$, $x \in C$. In
what follows, for simplicity of exposition, we shall consider the "strongly aperiodic case" $m=1$, that is
\begin{hyp}
\label{assum:smallset} There exist $\epsilon \in (0,1]$, a probability measure $\nu$  on $(\Xset,\Xsigma)$
such that $\nu(C) = 1$ and for all $x \in C$,
$A \in \Xsigma$, $P(x,A) \geq \epsilon \nu(A)$.
\end{hyp}
The general $m$ case can be straightforwardly, but to the price of heavy notations and calculus (considering for example easy extensions of i.i.d. theorem to the 1-dependent case), recovered from the proofs presented here. Fundamental to our methodology will be the regeneration technique (see \cite[chapter IV]{nummelin:1984}. The existence of small sets enables the use of the splitting
construction  to create atoms and to use regeneration methods, similar to those on countable spaces. In particular, each time the chain
reaches $C$, there is a possibility for the chain to regenerate. Each time the chain is at $x \in C$, a coin is tossed with probability of success $\epsilon$.
if the toss is successful, then the chain is moved according to the probability distribution $\nu$, otherwise, according to $(1-\epsilon)^{-1} \left\{ P(x,\cdot) - \epsilon \nu(\cdot) \right\}$.
Overall, the dynamic of the chain is not affected by this coin toss, but at each time the toss is successful, the chains regenerates with regeneration distribution $\nu$
independent from $x$. We denote by $\tau = \inf\{ k \geq 1, X_k \in C\}$ and $\sigma = \inf \{ k \geq 0, X_k \in C \}$ the first return and hitting
time to $C$ and by $\ctau= \inf \left\{ k \geq 1, (X_k,d_k) \in C \times \{1\} \right\}$ and $\csigma = \inf \left\{k \geq 0, (X_k,d_k) \in \Xset \times \{1\} \right\}$.
Let $f$ be a non-negative function and $r \in \Lambda$ a subgeometric sequence and $\mu$ a probability measure
on $(\Xset,\Xsigma)$. Our main result gives a bound to the $(f,r)$-modulated expectation of moments $\cPE_{\cmu} \left[ \sum_{k=1}^{\csigma} r(k) f(X_k) \right]$ of the regeneration time
(where $\cPE_{\cmu}$ is the expectation associated to the split chain; see below) in terms of the corresponding moment of
$\tPE_{\mu} \left[ \sum_{k=0}^{\tau} r(k) f(X_k) \right]$ and constants depending only and explicitly on $\epsilon$ and $\nu$ and on
the sequence $r$. Here, $\tPE_{\mu}$ denotes the
expectation associated to a Markov chain with initial distribution $\mu$ and moving according to $P$ outside $C$ and the residual kernel $(1-\epsilon)^{-1} \{ P(x,\cdot) - \epsilon \nu(\cdot) \}$
inside $C$.

Because finding bounds for $\tPE_{\mu} \left[ \sum_{k=0}^{\tau} r(k) f(X_k) \right]$  is not always easy, we will consider bounds for this
quantity derived from a "subgeometric" condition recently introduced in \cite{douc:fort:moulines:soulier:2004}, which might be seen, in the
subgeometrical case, as an analog to the Foster-Lyapunov drift condition for geometrically ergodic Markov Chains. We obtain, using these
drift conditions, explicit bounds for the $(f,r)$-modulated expectation of moments of the regeneration times in terms of the constants
in \textbf{A\ref{assum:smallset}}, the sequence $r$ and the constants appearing in the drift conditions.
With these results, we obtain  limit theorems for additive functionals and deviations inequalities, under conditions
which are easy to check.

\section{Bounds for regeneration time}
\label{sec:boundsregenerationtime}

We proceed by recalling the construction of the split chain \citep[Chapter 4]{nummelin:1984}.
For $x \in C$ and $A \in \Xsigma$ define the kernel $\RK$ as follows,
\begin{align}
\label{eq:residualkernel}
&\RK(x,A) = \begin{cases} (1-\epsilon \1_C(x) )^{-1} \left\{ P(x,A) - \epsilon \1_C(x) \nu(A) \right\} & 0 \leq \epsilon \1_C(x) < 1,\\
                                \delta_x(A) & \epsilon \1_C(x) =1
          \end{cases}
\end{align}
Define now, on the product space $\check{\Xset} = \Xset \times \{0,1\}$ equipped with the product $\sigma$-algebra $\Xsigma \otimes \mathcal{P}(0,1)$
where $\mathcal{P}(0,1) \eqdef \{ \emptyset, \{0\}, \{1\}, \{0,1\} \}$ the split kernel as follows:
\begin{align*}
& \cP(x,0; A \times \{0\}) = \int_A Q(x,dy) \{ 1 - \epsilon \1_C(y) \} && \cP(x,0; A \times \{1\}) = \epsilon Q(x,A \cap C)  \\
& \cP(x,1; A \times \{0\}) = \int_A \nu(dy) \{ 1 - \epsilon \1_C(y) \} &&  \cP(x,1; A \times \{1\}) = \epsilon \nu(A \cap C) \eqsp.
\end{align*}
For $\mu$ be a probability measure on $(\Xset,\Xsigma)$, define the split probability $\cmu$ on $(\Xset \times \{0,1\}, \Xsigma
\otimes \mathcal{P}(\{0,1\})$ by
\begin{align}
\label{eq:splitmeasure}
\cmu(A \times \{0\}) &= \int_A \{ 1-\epsilon \1_C(y) \} \mu(dy) \eqsp, && A \in \Xsigma \eqsp,\\
\cmu(A \times \{1\}) &= \epsilon \mu(A \cap C) \eqsp.
\end{align}
We denote by $\cPP_{\cmu}$ and $\cPE_{\cmu}$ the probability and the expectation on $ ( \Xset^\nset \times \{0,1\}^{\nset},$ $\Xsigma^\nset \otimes \mathcal{P}^{\otimes \nset}(\{0,1\}) )$
associated to the Markov chain $\{ X_n, d_n \}_{n \geq 0}$ with initial distribution $\cmu$ and transition kernel $\cP$. The definition of the split kernel
implies that
\begin{align}
\label{eq:bonnemarginale}
& \CPP[\cPP]{X_{n+1} \in A}{\mathcal{F}_n^X \vee \mathcal{F}_{n-1}^d} = P(X_n,A) \\
\label{eq:headprobability}
& \CPP[\cPP]{d_n=1}{\mathcal{F}_n^X \vee \mathcal{F}_{n-1}^d} = \epsilon \1_C(X_n) \\
\label{eq:onregenere}
& \CPP[\cPP]{X_{n+1} \in A}{\mathcal{F}_n^X \vee \mathcal{F}_{n-1}^d; d_n=1} = \nu(A) \eqsp,
\end{align}
where for $n \geq 0$, $\mathcal{F}_n^d = \sigma( d_k, k \leq n)$ and by convention $\mathcal{F}_{-1}^d = \{ \emptyset, \Omega)$.
Condition \eqref{eq:bonnemarginale} simply states that $\{X_n\}_{n \geq 0}$ is a Markov chain \wrt\ the filtration $(\mathcal{F}_n^X \vee \mathcal{F}_{n-1}^d, n \geq 0)$.
Condition \eqref{eq:headprobability} means that the probability of getting a head ($d_n=1$) as the $n$-th toss is equal to $\epsilon \1_C(X_n)$, independently of the
previous history $\mathcal{F}_{n-1}^X$ and of the $n-1$ previous toss.
Condition \eqref{eq:onregenere} says that, if head is obtained at the $n$-th toss ($d_n=1$), then the next transition obeys the transition law $\nu$ independently of
the past history of the chain and of the tosses. This means in particular that $\Xset \times \{1\}$ is a proper atom.
From conditions \eqref{eq:bonnemarginale}, \eqref{eq:headprobability} and \eqref{eq:onregenere}, we have
\[
\CPP[\cPP]{X_{n+1} \in A}{\mathcal{F}_n^X \vee \mathcal{F}_{n-1}^d; d_n=0}  = Q(X_n,A)  \eqsp.
\]
We denote respectively by $\tilde{\PP}_\mu$ and $\tPE_\mu$ the probability and the expectation
on $(\Xset^\nset, \Xsigma^{\otimes \nset})$ of a Markov chain with initial distribution $\mu$ and transition kernel $Q$.

Denote by $\{ \sigma_j \}_{j \geq 0}$ are the successive hitting times
of  $\{ X_{n} \}$ to the set $C$
\begin{equation}
\label{eq:successivehittingtime}
\sigma_0 \eqdef \inf \left\{ n \geq 0, X_{n} \in C \right \} \quad \text{and} \quad \sigma_j = \inf \{n > \sigma_{j-1}, X_{n} \in C \}, \quad j \geq 1 \eqsp,
\end{equation}
and by $N_n$  the number of visits of  $\{ X_{n}\}$ to the set $C$ before time $n$,
\begin{equation}
\label{eq:definitionNn}
N_n  = \sum_{i=0}^n \1_C(X_{n}) =\sum_{j=0}^\infty \1_{ \{ \sigma_j \leq  n \}}
\end{equation}
Define by $\csigma$ the hitting time of the atom of the split chain $\Xset \times \{1\}$,
\begin{equation}
\label{eq:definitioncsigma}
\csigma \eqdef \left\{ k \geq 0, d_k = 1 \right\} \eqsp.
\end{equation}
The stopping time $\csigma$ is a regeneration time and $\nu$ is a regeneration measure, \ie\ the distribution of $X_n$ conditional to $\csigma=n$ is
$\nu$ independently of the past history of the chain.
The following proposition relates the functionals of the regeneration time under the probability associated to the split chain $\cPP_{\cmu}$ to the corresponding functionals of the chain $\{ X_n \}$ under the probability $\tPP_{\mu}$.
\begin{proposition}
\label{prop:keyproposition}
Assume \textbf{A\ref{assum:smallset}}.
Let $\mu$ a probability measure on $(\Xset,\Xsigma)$. Let $\{ \xi_n \}$ be a non-negative $\mathcal{F}^X$-adapted process and let  $S$ be a
$\mathcal{F}^X$-stopping time. Then,
\begin{align}
\label{eq:keyrelation}
&\cPE_{\check{\mu}} \left[ \xi_{S} \1_{ \{ S < \csigma \} } \right] = \tPE_\mu \left[ \xi_{S} (1-\epsilon)^{N_S} \1_{ \{ S < \infty \} } \right]\eqsp, \\
\label{eq:boundreturntime}
& \cPE_{\cmu} \left[ \xi_{\csigma} \1_{ \{\csigma < \infty \}} \right] = \epsilon \sum_{j=0}^\infty (1-\epsilon)^j \tPE_{\mu} \left[\xi_{\sigma_{j}} \1_{ \{ \sigma_{j} < \infty \}} \right] \eqsp.
\end{align}
\end{proposition}
The proof is given in the Appendix \ref{sec:proof:keyproposition}.

We will now apply the proposition above to functionals of the form $\xi_n : = \sum_{k=0}^n r(k) g(X_k)$ where $g$ is a non-negative function and
$r \in \Lambda$ is a sequence, to relate the bounds of the $(g,r)$-modulated expectation of moments of regeneration time to the
$(f,r)$-modulated expectation of moments of the hitting time.
\begin{proposition} \label{prop:key2}
  Assume \textbf{A\ref{assum:smallset}}. Let $\{r(n)\}_{n \geq 0}$ be a sequence such that, for some $K$, $r(n+m) \leq K r(n) r(m)$, for
  all $(n,m) \in \nset \times \nset$. Let $g: \Xset \to [1,\infty)$ be a measurable function.
  For $x \in \Xset$, define
\begin{equation}
\label{eq:definitionWrg}
W_{r,g} (x) \eqdef \tPE_{x} \left[\sum_{k=1}^{\tau} r(k) g(X_{k}) \right],
\end{equation}
Then, for any $x \in \Xset$,
\begin{multline} \label{eq:relfund}
  \cPE_{\check{\delta}_x} \left [ \sum_{k=0}^{\csigma} r(k) g(X_{k})  \right] \leq \\
   r(0) g(x) + W_{r,g}(x)\1_{C^c}(x) +  \epsilon^{-1} (1-\epsilon) \, K \, \left( \sup_C W_{r,g} \right) \; \cPE_{\check{\delta}_x} [r(\csigma)]\eqsp.
\end{multline}
\end{proposition}
If $g \equiv 1$ and $r(n)= \beta^n$, this proposition may be seen as an extension of \cite[Theorem 2.1]{roberts:tweedie:1999},
which relates the generating function of the regeneration time to that of the hitting time to $C$.
Subgeometric sequences $r \in \Lambda_0$ also satisfies the inequality $r(n+m) \leq r(n) r(m)$. There is however a striking difference
with geometric sequence. Whereas for a geometric sequence $\liminf_{n \to \infty} r(n) / \sum_{k=0}^n r(k) > 0$, for subgeometric sequence we have on the
contrary $\limsup_{n \to \infty} r(n) / \sum_{k=0}^n r(k) = 0$. This implies that, whereas $\cPE_{\check{\delta}_x} \left [ \sum_{k=0}^{\csigma} r(k) g(X_{k})  \right]$ and
$\cPE_{\check{\delta}_x} [r(\csigma)]$ are of the same order of magnitude in the geometric case, the second is negligible compared to the
first one in the subgeometric case. In particular,
\begin{corollary} \label{coro:key2}
Assume \textbf{A\ref{assum:smallset}}. For any function $g: \Xset \to [0,\infty)$, there exists a constant $b_g$ (depending only and explicitly on $\epsilon$ and $\sup_C W_{1,g}$) such that
\begin{equation}
\cPE_{{\check{\delta}_x}} \left[ \sum_{k=0}^{\csigma} g(X_{k}) \right]\leq g(x) + W_{1,g}(x)\1_{C^c}(x)+ b_g \eqsp . \label{eq:requaltoone}
\end{equation}
For any $r\in\Lambda_0$ and $\delta > 0$, there exists a constant $b_r$ (depending only and explicitly on $\epsilon$, $\delta$, $r$ and
$\sup_C W_{r,1}$) such that
\begin{equation}
\cPE_{{\check{\delta}_x}} \left[ \sum_{k=0}^{\csigma} r(k) \right] \leq (1+\delta) W_{r,1}(x)\1_{C^c}(x)  + b_r \eqsp.\label{eq:relfund1}
\end{equation}
\end{corollary}
In general, of course, $\sup_C W_{1,g}$ and $\sup_C W_{r,1}$  is not easy to find analytically and, as in other approaches to this problem,
we will consider bounds on these quantities using "subgeometric drift" conditions as introduced in \cite{douc:fort:moulines:soulier:2004}, generalising a condition implying rieamnnian convergence stated in \cite{jarner:roberts:2001}
(see also \cite{fort:moulines:2000}). This condition may be seen  as an  analogue for
subgeometrically ergodic Markov chain of the Foster-Lyapunov condition for geometrically ergodic Markov chain.
\begin{hyp}
\label{assum:driftcondition} There exist a concave, non decreasing,
differentiable function $\varphi: [1,+\infty) \to \rset^+$,  a measurable function $V: \Xset \to [1,\infty)$ and
positive constants $b$ satisfying $\varphi(1) > 0$, $\lim_{v\to\infty} \varphi(v) = \infty$, $\lim_{v\to\infty}\varphi'(v) =
0$, $\sup_{x \in C} V(x) <\infty$ and
\[
P V \leq V - \varphi \circ V + b \1_C \eqsp, \\
\]
where the set $C$ is given in \textbf{A\ref{assum:smallset}}.
\end{hyp}
This drift condition has been checked in a large number of examples arising for example in queueing theory, Markov Chain Monte Carlo, time-series analysis (see for example
\cite{jarner:roberts:2001},\cite{douc:fort:moulines:soulier:2004}).
Examples of functions $\varphi$ satisfying \textbf{A\ref{assum:driftcondition}} include of course polynomial functions $\varphi(v) = (v+1)^\alpha$ for $\alpha \in (0,1)$ but also more general functions
like $\varphi(v) = \log^\alpha(v+1)$ for some $\alpha > 0$, or $\varphi(v) = (v+d) / \log(v+d)^\alpha$, for some $\alpha > 0$ and sufficiently large constant $d$.
We refer to \cite{douc:fort:moulines:soulier:2004} for precise statements giving both drift functions and rate $\varphi$ for these examples. Define
\begin{gather} \label{eq:defhphi}
\Hphi(v) \eqdef \int_1^v \frac{dx}{\varphi(x)}.
\end{gather}
The function $\Hphi : [1,\infty) \to [0,\infty)$ is  increasing and $\lim_{v \to \infty} \Hphi(v) = \infty$ (see \cite[Section 2]{douc:fort:moulines:soulier:2004}).
Define, for $u \in [1,\infty)$,
\begin{equation}
r_\varphi(u) \eqdef  \varphi \circ \Hphi^{-1}(u) / \varphi \circ \Hphi^{-1}(0), \label{eq:ratefunction}
\end{equation}
where $\Hphi^{-1}$ is the inverse of $\Hphi$. The function $u \mapsto r_\varphi(u)$ is log-concave and thus the sequence
$\{ r_\varphi(k) \}$ is subgeometric.
Polynomial functions $\varphi(v) = v^\alpha$, $\alpha \in (0,1)$ are associated to polynomial sequences $r_\varphi(k) =(1+(1-\alpha)k)^{\alpha/(1-\alpha)}$.
Functions like $\varphi(v) = c (v+d) /\log^\alpha(v+d)$ ($\alpha \in (0,1)$ and  sufficiently large $d$) are associated to subexponential sequences,
\[
r_\varphi(n) \asymp n^{-\alpha/(1+\alpha)} \exp\left( \{c(1+\alpha)
  n\}^{1/(1+\alpha)} \right) \eqsp.
\]
 where for two
sequences $\{u_k \}$ and $\{v_k\}$ of positive numbers, $u_k \asymp v_k$ if
$$
0 < \liminf_{k \to \infty} \frac{u_k}{v_k} \leq \limsup_{k \to \infty} \frac{u_k}{v_k} < \infty \eqsp.
$$
\cite[Proposition 2.2]{douc:fort:moulines:soulier:2004} shows that, under \textbf{A\ref{assum:smallset}-\ref{assum:driftcondition}},
for all $x \in \Xset$,
\begin{align}
\label{eq:boundmodulated-moment1}
& \PE_x \left [ \sum_{k=0}^{\tau_C-1} \varphi \circ V(X_k) \right ] \leq V(x) + b \1_C(x) \eqsp, \\
\label{eq:boundmodulated-moment2}
& \PE_x \left [ \sum_{k=0}^{\tau_C-1} r_\varphi(k) \right ] \leq \left\{ V(x) -1 +  b r_\varphi(1) \1_C(x) \right\} / \varphi(1) \eqsp.
\end{align}
This implies, using \cite{tuominen:tweedie:1994} that a Markov Chain satisfying \textbf{A\ref{assum:smallset}-\ref{assum:driftcondition}} is
both $(1,r_\varphi)$- and $(f,1)$-ergodic.
Denote by $\mathcal{G}(\varphi)$ the set of measurable functions
satisfying:
\begin{equation}
\label{eq:definitionGphi}
\mathcal{G}(\varphi) \eqdef \left\{ \psi: [1,\infty) \to \rset,  \text{$\psi$ is non decreasing, $\psi / \varphi$ is non increasing} \right\}\eqsp.
\end{equation}
Similarly to \eqref{eq:defhphi}, for all $\psi \in \mathcal{G}(\varphi)$, define the function
\begin{equation}
\label{eq:definitionLpsi}
\Hphi[\psi]: v \mapsto \int_{1}^v \frac{\psi}{\varphi}(u) du \eqsp.
\end{equation}
The function $\Hphi[\psi]$ is concave, non decreasing and, because
$[\psi/\varphi](u) \leq [\psi/\varphi](1)$, $\Hphi[\psi](u)\leq [\psi/\varphi](1) \, (u-1)$ for all $u \geq 1$.
The results of the previous section are used to derive explicit bounds  for
\begin{equation*}
   \cPE_{{\check{\delta}_x}} \left[ \sum_{k=0}^{\csigma} \psi \circ V(X_k) \right]\quad \mbox{and}\quad  \cPE_{{\check{\delta}_x}} \left[ \sum_{k=0}^{\csigma} r_\varphi(k) \right]
\end{equation*}
where $\psi$ is any function in $\mathcal{G}(\varphi)$.
The following theorem, proved in section \ref{sec:proof:theo:mainresult}, establishes bounds for the modulated moment of the excursion of the split chain to the atom $\Xset \times \{1\}$
as a function of the drift condition.
\begin{theorem}
\label{theo:mainresult} Assume
\textbf{A\ref{assum:smallset}-\ref{assum:driftcondition}}. Then,
there exists finite constant $B_\psi$ (depending only and explicitly on the constants appearing in the assumptions) such that for all $x \in \Xset$,
$\psi \in \mathcal{G}(\varphi)$,
\begin{equation}
\cPE_{{\check{\delta}_x}} \left[ \sum_{k=0}^{\csigma} \psi \circ V(X_k) \right]\leq \Hphi[\psi]\circ V(x) \1_{C^c}(x) + B_\psi  \eqsp, \label{eq:majosigcheckpsi}
\end{equation}
For any $\delta > 0$, there exists a finite constant $B_\varphi$ (depending only and explicitly on the constants appearing in the assumptions and $\delta > 0$)
such that
\begin{equation}
\cPE_{{\check{\delta}_x}} \left[ \sum_{k=0}^{\csigma} r_\varphi(k) \right] \leq (1 + \delta) V(x) \1_{C^c}(x)+ B_\varphi \label{eq:majosigcheckr}\eqsp,
\end{equation}
For any $c \in (0,1)$ and $K \geq 1$, there exists a finite constant $\kappa$
(depending only and explicitly on the constants appearing in the assumptions)
such that for any $\psi \in \GG(\varphi)$, and  $x \in \Xset$,
\begin{equation}
\label{eq:explicittail}
\cPP_{{\check{\delta}_x}} \left( \sum_{k=0}^{\csigma} \psi \circ V (X_{k}) \geq M \right) \leq \\ \kappa  \left[ \frac{1}{\Phi^{-1} \left\{ c M/ \psi(K) \right\}} + \frac{\Hphi[\psi](K)+1}{(1-c) M K}\right] V(x)\eqsp,
\end{equation}
\end{theorem}
The rates of convergence for the tail of the excursions may be obtained by optimizing the choice of the constant $K$ with respect to $M$.
As an illustration, consider first the case where $\psi \equiv 1$. Since $\lim_{s \to \infty} \varphi(s) = \infty$, then
$$
\lim_{K \to \infty} \frac{\Phi(K)}{K}= \lim_{K \to \infty} \frac{1}{K} \int_1^K \frac{ds}{\varphi(s)} = 0 \eqsp.
$$
Therefore, by letting $K \to \infty$ in the \rhs\ of \eqref{eq:explicittail} and then, taking $c=1$,
$$
\cPP_{{\check{\delta}_x}} \left(\csigma \geq M \right) \leq \kappa V(x) / \Phi^{-1}(M) \eqsp.
$$
Note that this bound could have been obtained directly by using the Markov inequality with the bound \eqref{eq:majosigcheckr} of the $f$-modulated moment of the
excursion. Consider now the
case: $\psi \equiv \varphi$. By construction, for any $K \geq 1$, $(\Hphi[\psi](K)+1)/K \leq 1$ and for any positive $u$, $\Phi^{-1}(u) \geq \varphi(1) u+1$.
Taking $K=1$ in \eqref{eq:explicittail}, Theorem \ref{theo:mainresult} shows that, for some constant $\kappa$,
$$
\cPP_{{\check{\delta}_x}} \left( \sum_{k=0}^{\csigma} \psi \circ V (X_{k}) \geq M \right) \leq \kappa V(x) /M \eqsp,
$$
which could have been again deduced from the Markov inequality applied to the bound for the excursion \eqref{eq:majosigcheckpsi}. The expression \eqref{eq:explicittail}
thus allows to retrieve these two extreme situations. Eq.~\eqref{eq:explicittail} also allows to interpolate the rates for functions growing more slowly than $\varphi \circ V$.

We give now two examples of convergence rates derived from the previous theorem by balancing  the two terms of the \rhs\ appearing in \eqref{theo:mainresult}.

\paragraph{{\bf Polynomial ergodicity.}} By Eqs (\ref{eq:defhphi}) and (\ref{eq:ratefunction}), if $\varphi(v)=  v^\alpha$ (with $\alpha \in (0,1)$), then
$r_\varphi(k)= \left( 1+(1-\alpha) k \right)^{\alpha/(1-\alpha)}$ and $\Phi^{-1}(u) \asymp (1-\alpha)^{1/(1-\alpha)} u^{1/(1-\alpha)}$ as $u \to \infty$.
Choose $\beta \in (0,\alpha)$ and set $\psi(u)=u^\beta$. Then, $\Hphi[\psi](v) = (1+\beta-\alpha)^{-1}(v^{1+\beta-\alpha}-1)$ and the optimal rate in the
\rhs\ of \eqref{eq:explicittail} is obtained by setting $K= M^{\frac{\alpha}{\beta + (\alpha-\beta)(1-\alpha)}}$. With this choice of $K$,
\eqref{eq:explicittail} implies that
$$
\cPP_{{\check{\delta}_x}} \left[ \sum_{k=0}^{\csigma} V^\beta(X_k) \geq M \right]\leq \kappa V(x) M^{-\frac{\alpha}{\beta + (\alpha-\beta)(1-\alpha)}} \eqsp.
$$
This bound shows how the rate of convergence of the tail depends on the tail behavior of the function $g$ and of the mixing rate of the Markov Chain.

\paragraph{\bf Subexponential ergodicity.} Assume that $\varphi(v)=c(v+d) (\log (v+d))^{-\alpha}$ for some positive constants $c$ and
$\alpha$ and sufficiently large $d$. Then  $\Phi^{-1}(k) \asymp  \rme^{(c(1+\alpha)k)^{1/(1+\alpha)}}$.
Choose for example $\psi(v)= |\log|^{\beta} (1+v)$, $v \in \rset^+$. By optimising the bound  \wrt\ $K$, \eqref{eq:explicittail} yields:
$$
 \cPP_{{\check{\delta}_x}} \left[ \sum_{k=0}^{\csigma} |\log|^\beta [V(X_{k})]  \geq M \right]\leq \kappa \rme^{-c M^{\frac{1}{1+\alpha+\beta}}} V(x)\eqsp,
$$
for some  constants $c$ and $C$ which does not depend of $\beta$ or $M$. Similarly, for $\psi(v)=(1+v)^\beta$ with $\beta \in (0,1)$, there
exists a constant $\kappa < \infty$,
$$
 \cPP_{{\check{\delta}_x}} \left[ \sum_{k=0}^{\csigma} V^\beta(X_{k})  \geq M \right]\leq \kappa M^{\frac{-1}{\beta}} \log^{\frac{2 \alpha \beta -1 + \beta - \alpha}{\beta}}(M) V(x)\eqsp.
$$

\section{Central Limit Theorem and Berry-Esseén Bounds}
\label{sec:CLT-BerryEsseen}
As a first elementary application of the results obtained in the previous section, we will derive conditions upon which a Central Limit Theorem holds for the normalized sum
$S_n(f) \eqdef $ $n^{-1/2} \sum_{i=1}^n (f(X_k) - \pi(f))$ where $\pi$ is the stationary distribution for the chain. For $u$, $v$ two vectors of $\rset^d$,
denote by $\pscal{u}{v}$ the standard scalar product and $\|u\|= \left( \pscal{u}{u}\right)^{1/2}$ the associated norm.

\begin{theorem}
\label{theo:CLTgen}
Assume \textbf{A\ref{assum:smallset}-\ref{assum:driftcondition}}. Let $\psi$ be a function such that
$\psi^2$ and $\psi \Hphi[\psi]$ belong to $\mathcal{G}(\varphi)$.
Then, for any function  $f: \Xset \to \rset$ such that $\sup_{\Xset} \frac{|f|}{\psi \circ V} < \infty$,
$$
\int f^2 d \pi  + \int | f | \sum_{k=1}^\infty |P^k f - \pi(f)| d \pi < \infty \eqsp.
$$
If in addition  $\stdlimit^2(f)>0$, where
\begin{equation}
\label{eq:defSigma}
\stdlimit^2(f) \eqdef \int \{f - \pi(f) \}^2   d \pi + 2  \int  f \sum_{k=1}^\infty P^k \{f - \pi(f)\} d \pi \eqsp,
\end{equation}
then, for any initial probability measure $\mu$ on $(\Xset,\Xsigma)$ satisfying $\mu(\Hphi[\psi]) < \infty$,
 $\sqrt{n} S_n(f)$ converges in distribution to a zero-mean Gaussian variable with variance $\stdlimit^2(f)$.
\end{theorem}

\paragraph{\textbf{Polynomial ergodicity}:}
Assume that $\varphi(v)=v^\alpha$ for some $\alpha \in (1/2,1)$ and choose $\psi(v)=v^\beta$ for some $\beta \in [0,\alpha]$.
Then,  $\Hphi[\psi](v)=(1+\beta-\alpha)^{-1}(v^{1+\beta-\alpha}-1)$ and the conditions of theorem \ref{theo:CLTgen}
are satisfied if $\alpha > 1/2$ and $\beta \in [0,\alpha-1/2]$. This is equivalent to the condition used in the CLT
\cite[Theorem 4.4]{jarner:roberts:2001} for polynomially ergodic Markov chains.
Note that, if $\alpha < 1/2$, then the moment of order two of the hitting time $\csigma$ is not necessarily finite,
and the CLT does not necessarily holds in this case.

\paragraph{\textbf{Subexponential ergodicity}:} Theorem \ref{theo:CLTgen} allows to derive a CLT under conditions which imply subexponential convergence.
Assume that $\varphi(v)=(d+v) \log^{-\alpha} (d+v)$, for some $\alpha>0$ and sufficiently large $d$. The condition of Theorem \ref{theo:CLTgen} are satisfied
for  $\psi(v)\propto v^{1/2} \{ \log (v) \}^{-(\alpha+\delta)}$ for $\delta>0$.

By strengthening the assumptions, it is possible to establish a Berry-Esseén Theorem with an explicit control of the constants.
\begin{theorem}
\label{theo:berry-esseen-bound}
In addition to the assumptions of Theorem \ref{theo:CLTgen}, suppose that the functions $\psi^3$, $\psi^2 \Hphi[\psi]$ and  $\psi \Hphi[{\psi \Hphi[\psi]}]$ belong
to $\mathcal{G}(\varphi)$. Let $\mu$ be a probability measure on $(\Xset,\Xsigma)$ such that $\mu(\Hphi[\psi]) < \infty$.
Then, there exist a constant $\kappa$ depending only and explicitly on the constants appearing in the assumptions (\textbf{A}\ref{assum:smallset}-\ref{assum:driftcondition})
and on the probability measure $\mu$ such that, for any function $f: \Xset \to \rset$ such that $\sup_{\Xset} \frac{|f|}{\psi \circ V} < \infty$ and $\stdlimit^2(f) > 0$,
\begin{equation}
\label{eq:Berry-Esseen-Bound}
\sup_{t} \left| \PP_{\mu} \left(  n^{-1/2}  S_n(f) / \stdlimit(f) \leq t\right)- G(t)  \right| \leq \kappa n^{-1/2} \eqsp,
\end{equation}
where $G$ is the standard normal distribution function.
\end{theorem}
Berry-Esseen theorems have been obtained for Harris-recurrent Markov chains under moment and strongly mixing conditions by \cite{bolthausen:1982}.
The use of the results obtained above allow to check these conditions directly from the drift condition. A side result, which is not fully
exploited here because of the lack of space, is the availability of an explicit computable expression for the constant $\kappa$,
which allows to investigate to assess deviation of the normalized sum for finite sample. 
This provides an other mean to get "honest" evaluation of the convergence of the Markov chain, under conditions which are less stringent than the ones outlined in \cite{jones:hobert:2001},
based on total variation distance.
It is interesting to compare our conditions with those derived in \cite[Theorem 1]{bolthausen:1982}, in the polynomial case, \ie\ $\varphi(v) =v^\alpha$
$\alpha \in (0,1)$. It is straightforward to verify that the conditions of  the Theorem \ref{theo:berry-esseen-bound} are satisfied by $\psi(v) = v^\beta$ if
$\alpha > 2/3$ and $\beta \in [0,\alpha-2/3]$. On the other hand, the strong mixing rate of this chain is $r(n) = n^{-\alpha/(1-\alpha)}$ (see
\cite{douc:fort:moulines:soulier:2004} and the maximum value of $p$ such that $\pi( V^{p\beta}) < \infty$ is $p= \alpha/\beta$. The Bolthausen condition
$\sum_{k=1}^\infty k^{(p+3)/(p-3)} r(n) < \infty$, is therefore satisfied again if $\alpha > 2/3$ and $\beta \in [0,\alpha-2/3)$, the value $\alpha-2/3$
being this time excluded.


\section{Moderate deviations}
\label{sec:MDP}

The main goal of this section is to generalize the MDP result of Djellout-Guillin \cite{djellout:guillin:2001} from the atomic case to the 1-small set case.
We will indicate in the proof the easy modifications needed to cover the general case.

\subsection{Moderate deviations for bounded functions}
We first consider MDP for bounded mapping, including non separable case (the functional empirical process and the trajectorial case).
\begin{theorem}
\label{theo:MDPgen}
Assume conditions \textbf{A\ref{assum:smallset}-\ref{assum:driftcondition}}.
Then, for all sequence $\{ b_n \}$ satisfying, for any $\varepsilon > 0$,
\begin{align}
&\lim_{n \to \infty} \left( \frac{\sqrt n}{b_n} + \frac{b_n}{n} \right) =0, \label{eq:speedcond0} \\
&\lim_{n \to \infty} \frac{n}{b_n^2} \log \left(\frac{n}{\Phi^{-1}(\varepsilon b_n)} \right) = - \infty \eqsp, \label{eq:speedcond1}
\end{align}
for all initial measure $\mu$ satisfying $\mu(V) < \infty$, for all bounded measurable function $f: \Xset \rightarrow \rset^d$ such that $\pi(f)=0$ and
for all closed set $F \subset \rset^d$, we have
\[
\limsup_{n\to\infty}{n \over b^2_n}\log \PP_\mu\left( \frac{1}{b_n} \sum_{k=0}^{n-1} f(X_k) \in F \right)
\le - \inf_{x\in F} J_f(x) \eqsp,
\]
where $J_f$ is a good rate function, defined by
\begin{equation}
\label{eq:defJf}
J_f(x) \eqdef \sup_{\lambda \in \rset^d} \left( \pscal{\lambda}{x} -(1/2)  \stdlimit^2(\lambda,f) \right),
\end{equation}
and $\stdlimit^2$ is defined by \eqref{eq:defSigma}.
\end{theorem}
The proof is given in section \ref{sec:proof:theo:MDPgen}.
\cite{acosta:1997} proved that the moderate deviation lower bound holds for all bounded function and all initial measure
provided that the chain is ergodic of degree 2, \ie\ for all set $B \in \Xsigma$ such that $\pi(B) > 0$,
$\int_B \PE_x[\tau_B^2] \pi(dx) < \infty$, where $\tau_B \eqdef \inf\{ k \geq 1, X_k \in B \}$ is the return-time to the set $B$.
It turns out that, under the assumptions \textbf{A\ref{assum:smallset}-\ref{assum:driftcondition}}, the  condition \eqref{eq:speedcond0}-\eqref{eq:speedcond1}
implies that the Markov chain is ergodic of degree 2. Note indeed that the conditions
\eqref{eq:speedcond0}-\eqref{eq:speedcond1} implies that $\lim_{k \to \infty} k / r_\varphi(k)= 0$. The definition \eqref{eq:ratefunction} of $\{ r_\varphi(k) \}$
implies that for some positive $c$, $\varphi(v) \geq c \sqrt{v}$, for any $v \in [1,\infty)$ and Lemma \ref{lem:characterization:ergodic:degree:two}
(stated and proved in section \ref{sec:proof:theo:MDPgen})
shows that this condition implies that the chain is ergodic of degree two. Thus, Theorem \ref{theo:MDPgen} together with
\cite[Theorem 3.1]{acosta:1997} establish the full MDP for  bounded additive functionals.

Condition \eqref{eq:speedcond0}-\eqref{eq:speedcond1} linking ergodicity and speed of the MDP may be seen as the counterpart for  Markov chains  of the condition of \cite{ledoux:1992} for the MDP of
\iid\ random variable linking the tail of this random variable with the speed of the MDP. Let us give examples of the range of speed of the MDP allowed as the function of the ergodicity rate.

\paragraph{\textbf{Polynomial ergodicity}:} By Eqs (\ref{eq:defhphi}) and (\ref{eq:ratefunction}), if $\varphi(v)=v^\alpha$ (with $\alpha \in (0,1)$), then
$r_\varphi(k) \asymp k^{\alpha/(1-\alpha)}$ and $\Phi^{-1}(k) \asymp k^{1/(1-\alpha)}$. Therefore, condition \eqref{eq:speedcond0}-\eqref{eq:speedcond1} is fulfilled as soon as
for any $\alpha \in (1/2,1)$ by any sequence $\{ b_n \}$ satisfying
$$\lim_{n \to \infty} \left\{ \frac{\sqrt{n}}{b_n} + \frac{\sqrt{n\log n}}{b_n} \right\}= 0 \eqsp.$$
\paragraph{\textbf{Subexponential ergodicity}:} Assume that $\varphi(v)=(v+d) (\log (v+d))^{-\alpha}$ for some $\alpha > 0$ and sufficiently large $d$.
Then, $\Phi^{-1}(k) \asymp e^{c k^{1/(1+\alpha)}}$ for some constant $c$. The condition \eqref{eq:speedcond0}-\eqref{eq:speedcond1} is fulfilled by any speed sequence $\{ b_n\}$ satisfying
$$\lim_{n \to \infty} \left\{ \frac{\sqrt{n}}{b_n} +  \frac{b_n}{n^{1+\alpha \over 1+2\alpha}} \right\}= 0 \eqsp.$$

The result can be extended to the empiral measure of a Markov chain. Assume that $\Xset$ is a Polish space and denote by
$\mathsf{M}(\Xset)$ the set of finite Borel signed measures on $\Xset$. Denote by $B(\Xset)$  the collection of bounded measurable
functions on $\Xset$. We equip $\mathsf{M}(\Xset)$ with the smallest topology such that the maps $\nu \mapsto \int_\Xset f d \nu$ are continuous for
each $f \in B(\Xset)$, commonly referred to as  the $\tau$-topology. The $\sigma$-algebra $\mathcal{M}(\Xset)$ on $\mathsf{M}(\Xset)$
is defined to be the smallest $\sigma$-algebra
such that for each $f \in B(\Xset)$, the map $\nu \mapsto f d \nu$ is measurable.
Define  the empirical measure $L_n$ as
$$
L_n={1\over b_n}\sum_{k=0}^{n-1}(\delta_{X_k}-\pi) \eqsp.
$$
For any $B \in \mathcal{M}(\Xset)$, we denote by $\mathrm{int}_\tau(B)$ and $\mathrm{clos}_\tau(B)$ the interior and the closure of the set $B$ \wrt\
the $\tau$-topology.
\begin{theorem}
\label{theo:empirMeasurMdp}
Under the assumptions of Theorem \ref{theo:MDPgen}, for every
probability measure $\mu \in \mathsf{M}(\Xset)$ satisfying $\mu(V)< \infty$, and any $B \in \mathcal{M}(\Xset)$
\begin{align*}
&\lim \sup_n \frac{n}{b_n^2} \log \PP_\mu (L_n \in B) \leq -\inf_{\gamma \in \ \mathrm{clos}_\tau (B)} I_0(\gamma) \eqsp, \\
&\lim \inf_n \frac{n}{b_n^2} \log \PP_\mu (L_n\in B) \geq -\inf_{\gamma \in \mathrm{int}_\tau(B)} I_0(\gamma)
\end{align*}
where for $\gamma \in \mathsf{\Xset}$, setting $\bar{f} = f - \pi(f)$,
\begin{gather}
I_0(\gamma)= \sup_{f \in B(\Xset)} \left[ \int f d \gamma - \frac{1}{2}\left( \int \bar{f}^2 d \gamma + 2 \int \bar{f} \sum_{k=1}^\infty P^k \bar{f} d \pi \right)\right] \eqsp.
\end{gather}
\end{theorem}
The proof can be directly adapted from the proof of \cite[Theorem 3.2]{acosta:1997} and is omitted for brevity.
An explicit expression of the good rate function can be found in \cite[Theorem 4.1]{acosta:1997}.
Other MDP principles (for instance, for the supremum of the empirical process) can be obtained, using the results obtained previously by \cite{djellout:guillin:2001}.
To save space, we do not pursue in this direction.

\subsection{Moderate deviations for unbounded functionals of
Markov chains}

We give here conditions allowing to consider unbounded functions. These conditions
make a trade-off  between the ergodicity of the Markov Chain, the range of speed for which a moderate deviation
principle may be established and the control of the tails of the functions.

\begin{theorem} \label{theo:mdpunbounded1}
Assume \textbf{A\ref{assum:smallset}-\ref{assum:driftcondition}} and
that there exist a function $\psi \in\mathcal{G}(\varphi)$ and a
sequence $\{ K_n \}$ such that $\lim_{n \to \infty} K_n = \infty$
and, for any positive $\varepsilon$,
\begin{align}
\label{eq:theo:mdpunbounded1}
&\lim_{n\to\infty}{n\over b^2_n}\log\left({n\over \Phi^{-1}(\varepsilon b_n/\psi(K_n))}\right)=-\infty \eqsp, \\
\label{eq:theo:mdpunbounded2}
&\lim_{n\to\infty}{n\over b^2_n}\log\left({n \Hphi[\psi](K_n) \over  \varepsilon b_n K_n}\right)=-\infty \eqsp.
\end{align}
Then, for any initial distribution $\mu$ satisfying $\mu(V)<\infty$ and any measurable function $f: \Xset \to \rset^d$ such that
$\sup_\Xset \| f \| / \psi \circ V$, the sequence $\{  \sigma_n^2(\lambda,f,\mu) \}$ where
$$
\sigma_n^2(\lambda,f,\mu) \eqdef \PE_\mu \left[ \left( \frac{1}{n} \sum_{k=0}^{n-1} \{ f(X_k) - \pi(f) \}\right)^2\right] \eqsp,
$$
has a limit $\sigma^2(\lambda,f)$ (which does not depend on $\mu$)
and $\PP_\mu \left[ L_n(f)\in\cdot \right]$ satisfies a moderate
deviation principle with speed $b_n^2/n$ and good rate function
$\mathcal{J}_f$,
$$
\mathcal{J}_f(x)=\sup_{\lambda \in \rset^d} \left[
\pscal{\lambda}{x} -(1/2) \sigma^2(\lambda,f) \right] \eqsp.
$$
Moreover, if $\psi^2 + \psi \Hphi[\psi] \in \mathcal{G}(\varphi)$, then
$\sigma^2(\lambda,f)=\stdlimit^2(\pscal{\lambda}{f})$ and $\mathcal{J}_f=J_f$.
\end{theorem}


\paragraph{\textbf{Polynomial ergodicity}:} By Eqs (\ref{eq:defhphi}) and (\ref{eq:ratefunction}), if $\varphi(v)=v^\alpha$ (with $\alpha \in (1/2,1)$), then
$r_\varphi(k) \asymp k^{\alpha/(1-\alpha)}$ and $\Phi^{-1}(k) \asymp k^{1/(1-\alpha)}$. Choose $\psi(v)=v^\beta$ with $\beta<\alpha-1/2$. Then the MDP holds for
for any sequence $\{ b_n \}$ such that $\lim_{n \to \infty} \{ \frac{\sqrt{n}}{b_n} +  \frac{b_n}{\sqrt{n\log n}} \}= 0$. It is worthwhile to note that the speed which can be achieved
are the same than in the bounded case.

\paragraph{\textbf{Subexponential ergodicity}:} Assume now that $\varphi(v)=(v+d) (\log (v+d))^{-\alpha}$ for some $\alpha > 0$ and sufficiently marge $d$.  Then
Letting  $\psi(v)=(\log (1+v))^\beta$ for some $\beta>0$, then Theorem \ref{theo:mdpunbounded1} shows the MDP with speed $b_n=n^a$ for $a$ such that
$${1\over 2}<a<{\beta+1+\alpha\over 2\beta+1+2\alpha}.$$
Letting $\psi(v)=(1+v)^\beta$ with $\beta<1/2$, then Theorem \ref{theo:mdpunbounded1}
shows that the MDP principle holds for any sequence $\{b_n\}$ such that $\lim_{n \to \infty} \{ \frac{\sqrt{n}}{b_n} +  \frac{b_n}{\sqrt{n\log n}} \}= 0$.


\section{Deviation inequalities}
\label{sec:Deviation-Inequalities}
We now investigate some exponential deviation inequalities for
$\PP_\mu(\sum_{k=0}^n f(X_i)>\varepsilon n)$ valid for each $n$ where $f$ is a bounded and centered function w.r.t. $\pi$.
This is to be compared to Bernstein's inequality for i.i.d. variables or more precisely to the \cite{fuk:nagaev:1971} inequality  adapted to Markov chains, (as done in a previous work of
\cite{clemencon:2001}) except that in this paper, the Markov chain is not geometrically but sub-geometrically ergodic. Extensions to the case of unbounded functions
can be tackled using result of Theorem \ref{theo:mainresult}.

\begin{theorem}\label{theo:dev}
Assume that $f$ is bounded and centered with respect to $\pi$ and the assumptions of Theorem 1. Then, for
any initial measure $\mu$ satisfying $\mu(V) < \infty$, for any positive
$\varepsilon>0$ and $n>n_0(\varepsilon)$, there exists
$L,K$ (independent of $n$ and $\epsilon$) such that, for all positive $y$

$$\PP_\mu\left(\left\|\sum_{k=0}^{n-1}f(X_k)\right\|>\varepsilon n\right)\le \frac{Ln}{\Phi^{-1}\left(\frac{\varepsilon n}{\| f\|_\infty}\right)}+\frac{Ln}{\Phi^{-1}\left(\frac{y}{\| f\|_\infty}\right)}+e^{-\frac{n\varepsilon^2}{K\|f\|^2_\infty+\varepsilon y}}.$$

\end{theorem}
The proof is given in section \ref{sec:mdpdev}. Let us give a few comments on the obtained rate in some examples: with $\|f\|_\infty\le1$, for $n\ge n_1$
 
\begin{enumerate}
\item $\varphi(v)=(1+v)^{\alpha}$ for $\alpha\in(1/2,1)$, then there exists $K$
$$\PP_\mu\left(\left\|\sum_{k=0}^{n-1}f(X_k)\right\|>\varepsilon n\right)\le K\frac{\log(n)^{1\over 1-\alpha}}{\varepsilon^{2\over 1-\alpha}n^{\alpha\over 1-\alpha}}$$
\item $\varphi(v)=(1+v)\log(c+v)^{-\alpha}$ for positive $\alpha$, then there exists $K,L$
$$\PP_\mu\left(\left\|\sum_{k=0}^{n-1}f(X_k)\right\|>\varepsilon n\right)\le K~e^{-L(n\varepsilon)^{1\over 2+\alpha}}.$$
\end{enumerate}
The polynomial rate shown in the first case is better than the one derived by Rosenthal's inequality, and considering that we in fact only consider integrability assumptions, are not so far from optimal when considering stronger assumptions as weak Poincare inequalities. The subgeometric case is less satisfactory in the sense that when $\alpha$ is near 0, we hope to achieve a $n$ in the exponential (obtained for example via Cramer argument) whereas we obtained instead $\sqrt{n}$. The gap here, due to Fuk-Nagaev's inequality, is fullfilled only asymptotically via the moderate deviations result, and is left for deviation inequalities for further study.

\appendix

\section{Proof of  Propositions \ref{prop:keyproposition} and \ref{prop:key2}}
\label{sec:proof:keyproposition}
\newcommand{\condexp}[4]{#1_{#2}\left [ \left. #3 \right | #4 \right]}
\begin{proof}[Proof of the Proposition \ref{prop:keyproposition}]
We first prove by induction that for all $n \geq 0$ and all functions $f_0, \dots, f_n \in \pfunc{\Xset}$,
\begin{equation}
\label{eq:recurrence}
\cPE_{\cmu} \left[ \prod_{i=0}^n f_i(X_i) \1_{ \{ \csigma > n \}} \right] = \tPE_\mu \left[ \prod_{i=0}^n f_i(X_i) (1-\epsilon)^{N_n} \right]\eqsp.
\end{equation}
We first establish the result for $n=0$. For $f \in \pfunc{\Xset}$ we have
\begin{multline*}
\cPE_{\cmu} [ f(X_{0}) \1_{ \{ \csigma > 0 \}}] = \cPE_{\cmu}[  f(X_{0}) \1_{ \{d_0 = 0\}} ] = \\
(1-\epsilon) \int_C f(x) \mu(dx) + \int_{C^c} f(x) \mu(dx)= \int_\Xset \{ 1-\epsilon \1_C(x) \} f(x) \mu(dx),
\end{multline*}
Assume now that the result holds up to order $n$, for some $n \geq 0$. Similarly, for any $f \in \pfunc{\Xset}$,
\begin{align*}
&\cPE[f(X_{n+1} ) \1_{\{d_{n+1}=0\}} \mid \mathcal{F}_{n}^X \vee \mathcal{F}_n^d] \1_{ \{ d_n= 0\}}\\
&\quad = \cPE[ f(X_{n+1} ) \{ 1-\epsilon \1_C( X_{n+1}) \} \mid \mathcal{F}_{n}^X \vee \mathcal{F}_n^d] \1_{\{ d_n = 0\}} \\
&\quad = \tPE [ f(X_{n+1} ) \{ 1-\epsilon \1_C( X_{n+1}) \} \mid X_{n}]\1_{ \{ d_n= 0\}}
\end{align*}
Therefore, by the recurrence assumption,
\begin{align*}
&\cPE_{\cmu} \left[f_{n+1}(X_{n+1} )  \prod_{i=0}^n f_i(X_i) \1_{ \{ \csigma > n+1 \}} \right] \\
&\quad = \cPE_{\cmu} \left[ \tPE [f_{n+1}(X_{n+1} ) \{1-\epsilon \1_C(X_{n+1})\} \mid X_{n}] \; \prod_{i=0}^n f_i(X_i) \1_{ \{ \csigma > n \}} \right] \\
&\quad = \tPE_\mu \left[ \tPE[f_{n+1}(X_{n+1} ) \{ 1-\epsilon \1_C(X_{n+1}) \} \mid X_{n}] \; \prod_{i=0}^n f_i(X_i) (1-\epsilon)^{N_n} \right]\\
&\quad = \tPE_\mu \left[ f_{n+1}(X_{n+1}) \prod_{i=0}^n f_i(X_i) (1-\epsilon)^{N_{n+1}} \right]\eqsp,
\end{align*}
showing \eqref{eq:recurrence}. Therefore, the two measures on $(\Xset^{n+1}, \Xsigma^{\otimes (n+1)})$ defined respectively by
\begin{align*}
&A \mapsto \cPE_{\cmu} \left[ \1_A(X_0, \dots, X_n) \1_{ \{ \csigma \geq n  \}} \right] \quad \text{and} \\
&A \mapsto \tPE_\mu \left[ \1_A(X_0, \dots, X_1) (1-\epsilon)^{N_n} \right]
\end{align*}
are equal on the monotone class $\mathcal{C} \eqdef \{ A, A = A_0 \times \dots \times A_n, A_i \in \Xsigma \}$ for any $n$, and thus these two measures coincide on the
product $\sigma$-algebra.
The proof of \eqref{eq:keyrelation} follows upon conditioning upon the events $\{ S= n \}$. We now prove \eqref{eq:boundreturntime}.
By definition of the hitting time $\csigma$ to the atom $\Xset \times \{1 \}$,  $\xi_{\csigma} \1_{ \{\csigma < \infty \}}$ may be expressed as
$$
\xi_{\csigma} \1_{ \{\csigma < \infty\}} = \xi_{\sigma_0} \1_{ \{ d_{\sigma_0}= 1 \}} \1_{\{d_{\sigma_0} < \infty \}} +
\sum_{j=1}^\infty \xi_{\sigma_j} \1_{ \{ d_{\sigma_j} = 1 \}} \1_{ \{\sigma_{j-1} < \csigma \}} \1_{ \{ \sigma_j < \infty \}}.
$$
Note that
\[
\CPE[\cPE]{ \1_{\{ d_{\sigma_j}=1\}} \xi_{\sigma_j}}{\mathcal{F}^X_{\sigma_j}} \1_{\{\sigma_j < \infty\}}=
\epsilon \CPE[\cPE]{ \xi_{\sigma_j}}{\mathcal{F}^X_{\sigma_j} } \1_{\{ \sigma_j < \infty \}}\]
and $(1-\epsilon)^{N_{\sigma_j}} \1_{\{ \sigma_j < \infty \}}= (1-\epsilon)^{j+1} \1_{\{\sigma_j < \infty \}}$.
The proof follows from the identity
\begin{multline*}
\cPE_{\cmu}[ \xi_{\sigma_j} \1_{ \{ \sigma_j < \infty \}} \1_{ \{\sigma_{j-1} < \csigma \}}  ] =
\cPE_{\cmu}[ \cPE[ \xi_{\sigma_j} \1_{ \{ \sigma_j < \infty \}} \mid \mathcal{F}^X_{\sigma_{j-1}} \vee \mathcal{F}^d_{\sigma_{j-1}}] \1_{\{ d_{\sigma_{j-1}} = 0 \}} \1_{\{ \sigma_{j-1} < \csigma \}}] = \\
\tPE_\mu[ \tPE[ \xi_{\sigma_j} \1_{ \{ \sigma_j < \infty \}} \mid \mathcal{F}^X_{\sigma_{j-1}}] (1-\epsilon)^{N_{\sigma_{j-1}}} \1_{\{\sigma_{j-1} < \infty \}} ] =
(1-\epsilon)^j \tPE_\mu[ \xi_{\sigma_j} \1_{ \{ \sigma_j < \infty \}} ].
\end{multline*}
\end{proof}

\begin{proof}[Proof of the Proposition \ref{prop:key2}]
Without loss of generality we assume that $\sup_C W_{r,g} < \infty$ (otherwise the inequality is trivial).
  In the case $\epsilon=1$, Proposition \ref{prop:key2} is elementary
  since by Proposition \ref{prop:keyproposition}, it then holds that
$$
\cPE_{{\check{\delta}_x}} \left [\sum_{k=1}^{\csigma} r(k) {g}(X_{k}) \right]=W_{r,g}(x)\1_{C^c}(x).
$$
Consider now the case $\epsilon \in (0,1)$.  By applying Proposition \ref{prop:keyproposition}, we
obtain:
\begin{multline} \label{eq:applilemdouc}
  \cPE_{\check{\delta}_x} \left [ \sum_{k=1}^{\csigma} r(k) g(X_{k}) \right] =\\
  \epsilon W_{r,g}(x) \1_{C^c}(x)+ \epsilon \, \sum_{j=1}^\infty (1-\epsilon)^j
  \tPE_x \left[ \sum_{k=1}^{\sigma_j} r(k) {g}(X_{k}) \right].
\end{multline}
For $j \geq 1$, write
\begin{equation*}
  \tPE_x \left[ \sum_{k=1}^{\sigma_j} r(k) g(X_{k}) \right] =  W_{r,g} (x)\1_{C^c}(x)
  + \sum_{\ell=0}^{j-1} \tPE_x \left[\sum_{k=\sigma_\ell+1}^{\sigma_{\ell+1}} r(k) g(X_{k}) \right] \eqsp.
\end{equation*}
Under the stated assumptions, for all $n,m \geq 0$, $r(n+m) \leq K r(n)r(m)$.
This and the strong Markov property imply, for $x \in \{ W_{r,g} < \infty \}$~:
\begin{multline*}
\tPE_x \left[ \sum_{k=\sigma_\ell+1}^{\sigma_{\ell+1}} r(k) g(X_{k}) \right]  =
\tPE_x \left[ \sum_{k=1}^{\tau\circ\theta^{\sigma_\ell} } r(k+\sigma_\ell)    g(X_{k+\sigma_\ell}) \right] \\
\leq K  \tPE_x \left[ r(\sigma_\ell) W_{r,g} (X_{\sigma_\ell}) \right]  \leq  K \left(\sup_C W_{r,g} \right) \, \tPE_x[r(\sigma_\ell)],
\end{multline*}
where $\theta$ is the shift operator.
Plugging this bound into \eqref{eq:applilemdouc} and using again Proposition \ref{prop:keyproposition}, we obtain,
\begin{align*}
\cPE_{\check{\delta}_x} \left [ \sum_{k=1}^{\csigma} r(k) g(X_{k}) \right] & \leq W_{r,g}(x)\1_{C^c}(x) +  K \left( \sup_C W_{r,g} \right) \sum_{j=1}^\infty \epsilon (1-\epsilon)^j \sum_{\ell=0}^{j-1}  \tPE_x[r(\sigma_\ell)] \\
&  =  W_{r,g}(x)\1_{C^c}(x) +  K \left( \sup_C W_{r,g} \right) \sum_{\ell=0}^\infty (1-\epsilon)^{\ell+1} \, \tPE_x[r(\sigma_\ell)] \\
&  =  W_{r,g}(x)\1_{C^c}(x) + \epsilon^{-1} (1-\epsilon) K \, \left( \sup_C W_{r,g} \right) \cPE_{{\check{\delta}_x}}[r(\csigma)].
\end{align*}
\end{proof}

\begin{proof}[Proof of corollary \ref{coro:key2}]
For any $r \in \Lambda$, $\lim_{n \rightarrow \infty} r(n) / \sum_{k=1}^{n} r(k) = 0$.
As a consequence, for any $r(n) \in \Lambda$ and any $\delta > 0$, $N_{r,\delta}$ defined by
\begin{equation} \label{eq:defNr}
  N_{r,\delta} \eqdef \sup \left\{ n \geq 1, r(n) / \sum_{k=1}^{n}r(k) \geq \delta \right \} \eqsp,
\end{equation}
is finite. For any $n \geq 0$, the definition of $N_{r,\delta}$ implies
$r(n) \leq \delta  \sum_{k=1}^{n}   r(k) + r( N_{r,\delta}  ) \eqsp$.
Hence, for any $x \in \Xset$,
\[
\cPE_{{\check{\delta}_x}}[r(\csigma)]  \leq \delta \cPE_{{\check{\delta}_x}} \left [\sum_{k=1}^{\csigma} r(k)  \right] + r( N_{r,\delta} )
\]
The proof of \eqref{eq:relfund1} then follows by choosing $\delta$ sufficiently small so that $(1- \epsilon^{-1} (1-\epsilon) \sup_C W_{r,1} \delta)^{-1} \leq 1 + \delta$.
\end{proof}

\section{Proof of  Theorem \ref{theo:mainresult}}
\label{sec:proof:theo:mainresult}
We preface the proof by the following elementary lemma.
\begin{lemma}
\label{lem:driftpsi}
Assume \textbf{A}\ref{assum:driftcondition}. Then, for any $\psi \in \mathcal{G}(\varphi)$ there exists $b_\psi$
(depending only and explicitly on $b$, $\psi$ and $\varphi$) such that, for all $x \in \Xset$,
\begin{equation}
Q (x,\Hphi[\psi]\circ V) \leq \Hphi[\psi] \circ V(x) -\psi \circ V(x) + b_{\psi} \1_C(x) \label{eq:mdriftpsi}
\end{equation}
\end{lemma}

\begin{proof}
Since $\Hphi[\psi]$ is concave, differentiable, non decreasing, the Jensen inequality implies, for $x \not \in C$
\begin{multline*}
P [\Hphi[\psi]\circ V] \leq \Hphi[\psi](P V)\leq  \Hphi[\psi](V -\varphi \circ V) \\
\leq \Hphi[\psi] \circ V + \Hphi[\psi]'(V)(-\varphi \circ V) \leq \Hphi[\psi] \circ V -\psi \circ V
\end{multline*}
and
\[
\sup_C Q(x, \Hphi[\psi] \circ V) \leq \Hphi[\psi] \left[ (1-\epsilon)^{-1} \left\{ \sup_C PV - \epsilon \nu(V) \right\}\right] \eqsp.
\]
The proof follows.
\end{proof}

\begin{proof}
By Corollary \ref{coro:key2}, we may write
  \begin{equation}
  \cPE_{{\check{\delta}_x}} \left[ \sum_{k=0}^{\csigma} \psi \circ V(X_k) \right] \leq
  \tPE_{x} \left[ \sum_{k=0}^{\tau-1} \psi \circ V (X_{k}) \right] \1_{C^c}(x)  + \sup_C \psi \circ V + b_g \eqsp.
  \label{eq:firstcheck}
  \end{equation}
On the other hand, the comparison Theorem \cite[Theorem 11.3.1]{meyn:tweedie:1993} and the drift condition \eqref{eq:mdriftpsi} implies
that
\[
\tPE_{x} \left[ \sum_{k=0}^{\tau-1} \psi \circ V (X_{k}) \right] \1_{C^c}(x) \leq  \Hphi[\psi] \circ V(x) \1_{C^c}(x) \eqsp.
\]
The proof of \eqref{eq:majosigcheckpsi} follows. The proof of \eqref{eq:majosigcheckr} is along the same lines using
\eqref{eq:relfund1} instead of \eqref{eq:requaltoone}.

We now consider \eqref{eq:explicittail}. Define $\eta \eqdef \inf \{ k \geq 0, V(X_{k}) \geq K \}$.
We consider first the event $\left\{ \sum_{k=0}^{\csigma} \psi \circ V(X_k) \geq M, \eta \geq \csigma \right\}$,
on which $\psi \circ V(X_k)$ remains bounded by $\psi(K)$.
Therefore, on $\{ \eta \geq \csigma \}$, $\sum_{k=0}^{\csigma} \psi \circ V(X_k) \leq (\csigma + 1) \psi(K)$, which implies that
$$
\left\{ \sum\nolimits_{k=0}^{\csigma} \psi \circ V(X_{k}) \geq M, \eta \geq \csigma \right\} \subset \{  \csigma  \geq M / \psi(K) \} \eqsp.
$$
We now consider the complementary event: $\left\{ \sum_{k=0}^{\csigma} \psi \circ V(X_k) \geq M, \eta < \csigma \right\}$. We take $c \in (0,1)$,
Note that, if $\csigma < c M / \psi(K)$, then, $\sum_{k=0}^{\eta-1} \psi \circ V(X_k) \leq \eta \psi(K) \leq c M$ which implies that
$\sum_{k=\eta}^{\csigma} \psi \circ V(X_k) \geq (1-c)M$. Therefore,
\begin{multline*}
\left\{ \sum\nolimits_{k=0}^{\csigma} \psi \circ V(X_{k}) \geq M, \eta < \csigma  \right\} \subset \left\{  \csigma  \geq c M/ \psi(K) \right\} \\
\cup \left \{ \eta \leq  \csigma \leq c M / \psi(K) \eqsp, \sum\nolimits_{k=\eta}^{\csigma} \psi \circ V(X_k) \geq (1-c) M \right\} \eqsp.
\end{multline*}
Therefore,
\begin{multline}
\cPP_{{\check{\delta}_x}} \left( \sum\nolimits_{k=0}^{\csigma} \psi \circ V(X_{ k}) \geq M \right) \leq  2\cPP_{{\check{\delta}_x}} \left( \csigma \geq  c M/\psi(K)\right) + \\
\cPP_{{\check{\delta}_x}} \left ( \eta \leq \csigma  \leq c M / \psi(K) \eqsp , \sum\nolimits_{k=\eta}^{\csigma} \psi \circ V(X_{k}) \geq (1-c) M \right) \eqsp.
\label{eq:boundExcursion}
\end{multline}
The first term of the \rhs\ of (\ref{eq:boundExcursion}) is bounded using the Markov inequality with \eqref{eq:majosigcheckr},
$$
\cPP_{{\check{\delta}_x}} \left\{ \csigma \geq  c M / \psi(K) \right\} \leq
\frac{\cPE_{{\check{\delta}_x}} \left\{ \sum_{k=0}^{\csigma} r_\varphi(k) \right\}}{\Phi^{-1}\left\{ c  M / \psi(K) \right\}} \leq
\kappa_0  \frac{V(x)\1_{C^c}(x)+1}{\Phi^{-1}\left\{ c M/\psi(K) \right\}} \eqsp,
$$
for some finite constant $\kappa_0$. Similarly, the Markov inequality and the strong Markov property imply, using Eq.~\eqref{eq:majosigcheckpsi},
\begin{align*}
& \cPP_{{\check{\delta}_x}} \left\{\eta \leq  \csigma \leq cM /\psi(K), \sum_{k=\eta}^{\csigma} \psi \circ V(X_{k}) \geq (1-c)M  \right\}\\
& \quad \leq \frac{1}{(1-c)M}\cPE_{{\check{\delta}_x}}\left(\1_{\{ \eta \leq \csigma \}} \cPE\left\{ \left. \sum_{k=\eta}^{\csigma} \psi \circ V(X_k)\right| \mathcal{F}^X_\eta\right\}\right) \\
&\quad \leq \frac{\kappa_1}{(1-c)M} \cPE_{{\check{\delta}_x}}\left[ (\Hphi[\psi]  \circ V(X_{\eta})+1) \1_{\{\csigma  \geq \eta \}} \right] \eqsp,
\end{align*}
for some constant $\kappa_1$.
The function $u \mapsto \Hphi[\psi](u)/u$ is non-increasing. Therefore,
$(\Hphi[\psi] \circ V(X_{\eta})+1) \1_{\{ \eta < \infty \}} \leq $ $ K^{-1}(\Hphi[\psi](K)+1)V(X_{\eta}) \1_{\{\eta  < \infty\}}$,
which implies that
$$
\cPE_{{\check{\delta}_x}}\left[ (\Hphi[\psi]  \circ V(X_{\eta})+1) \1_{\{\csigma  \geq \eta \}} \right] \leq
\frac{(\Hphi[\psi](K)+1)}{K} \cPE_{{\check{\delta}_x}}\left[ V(X_{\eta}) \1_{\{\csigma \geq \eta \}} \right] \eqsp.
$$
We now prove that there exists a constant $\kappa_2$ such that, for any $x \in \Xset$,
\begin{equation}
\label{eq:keyrelation2}
\cPE_{{\check{\delta}_x}}\left[ V(X_{\eta})\1_{\{\eta\leq \csigma \}} \right] \leq \kappa_2 V(x) \eqsp.
\end{equation}
Since $\eta$ is $\mathcal{F}^X$-stopping time, using  proposition \ref{prop:keyproposition}, \eqref{eq:keyrelation}, we may write
\[
\cPE_{{\check{\delta}_x}}\left[ V(X_{\eta})\1_{\{\eta < \csigma \}} \right]=
\tPE_x \left[ V(X_\eta)  (1-\epsilon)^{N_\eta} \1_{\{\eta < \infty \}} \right] \eqsp.
\]
By conditioning upon the successive visit to the set $C$, the RHS of the previous equation may be expressed as
\begin{multline}
\label{eq:relation1}
\tPE_x \left[ V(X_\eta) (1-\epsilon)^{N_\eta} \1_{\{\eta < \infty \}} \right] =  \\
\tPE_x \left[ V(X_\eta) \1_{\{\eta < \sigma_0 \}} \right] + \sum_{j=1}^\infty (1-\epsilon)^{j} \tPE \left[ V(X_\eta) \1_{\{\sigma_{j-1} \leq \eta < \sigma_j\}} \right] \eqsp.
\end{multline}
Because $V(X_{\eta}) \1_{\{ \eta < \sigma_0 \}} \leq V(X_{\eta \wedge \sigma_0})$ and $\eta \wedge \sigma_0$ is a $\mathcal{F}^X$-stopping time,
the comparison Theorem (\cite[Theorem 11.3.1]{meyn:tweedie:1993}) implies that, under \textbf{A\ref{assum:driftcondition}},
\begin{equation}
\label{eq:relation2}
\tPE_x \left[ V(X_\eta) \1_{\{\eta < \sigma_0 \}} \right] \leq V(x) + b \1_C(x) \eqsp.
\end{equation}
Similarly, for any $j \geq 1$, we may write
\[
V(X_{\eta}) \1_{\{ \sigma_{j-1} \leq \eta < \sigma_j \}} \leq V(X_{\sigma_j \wedge \eta}) \1_{\{\sigma_{j-1} \leq \eta\}}
\leq V(X_{\tau \wedge \eta}) \circ \theta^{\sigma_{j-1}} \1_{\{ \sigma_{j-1} \leq \eta \}} \eqsp,
\]
and the comparison Theorem and the strong Markov property imply that
\begin{equation}
\label{eq:relation3}
\tPE_x \left[ V(X_{\eta}) \1_{\{ \sigma_{j-1} \leq \eta < \sigma_j \}} \right] \leq \left( \sup_C V + b \right) \eqsp.
\end{equation}
By combining the relations \eqref{eq:relation1}, \eqref{eq:relation2} and \eqref{eq:relation3}, we therefore obtain the bound
\[
\tPE_x \left[ V(X_\eta) (1-\epsilon)^{N_\eta} \1_{\{\eta < \infty \}} \right] \leq V(x) + b \1_C(x) + \frac{(1-\epsilon)}{\epsilon} \left\{ \sup_C V + b \right\} \eqsp,
\]
showing \eqref{eq:keyrelation2} and concluding the proof.
\end{proof}

\section{Proof of Theorem \ref{theo:CLTgen}}
\label{sec:proof:theo:CLTgen}
\begin{proof}[Proof of Theorem \ref{theo:CLTgen}]

By \cite[Theorem 17.3.6]{meyn:tweedie:1993}, we only need to check that $I= \cPE_{\check{\nu}} \left\{ \left(\sum_{k=0}^{\csigma} |f|(X_k)\right)^2 \right\}< \infty$.
We may write $I = I_1 + 2 I_2$ where the two terms $I_1$ and $I_2$ are respectively defined by
\begin{align*}
I_1 &\eqdef \cPE_{\check{\nu}}\left[\sum_{k=0}^{\csigma} f^2(X_k) \right]\\
I_2 &\eqdef \cPE_{\check{\nu}}\left[\sum_{k=0}^{\csigma} |f|(X_{k})  \sum_{\ell=k+1}^{\csigma} |f|(X_{\ell}) \right]  \\
    &=\cPE_{\check{\nu}}\left[\sum_{k=0}^{\csigma} |f|(X_{k}) \cPE_{X_{k},d_k} \left\{ \sum_{\ell=0}^{\csigma} |f|(X_{\ell})\right\} \right ]
\end{align*}
The proof follows using Theorem \ref{theo:mainresult}.
\end{proof}

\section{Proof of Theorem \ref{theo:MDPgen}}
\label{sec:proof:theo:MDPgen}
\begin{lemma}
\label{lem:characterization:ergodic:degree:two}
Assume that \textbf{A\ref{assum:smallset}-\ref{assum:driftcondition}} hold for some function $\varphi$ such that $\inf_{v \in [1,\infty)}$ $\frac{\varphi(v)}{\sqrt{v}} > 0$. Then, the chain
is ergodic of degree two.
\end{lemma}
\begin{proof}
Recall that for a phi-irreducible Markov Chain, the stationary distribution $\pi$ is a maximal irreducibility measure (see for instance \cite[Proposition 10.4.9]{meyn:tweedie:1993}),
Therefore any set $C \in \Xsigma$  such that $\pi(B) > 0$ is accessible. In addition, for any non-negative measurable function $f$,
$\pi(f) = $ $\int_B \pi(dx) \PE_x\left( \sum_{k=0}^{\tau_B-1} f(X_k) \right)$. A direct calculation shows that
$$
\PE_x[\tau_B^2] = 2 \PE_x\left[ \sum_{k=0}^{\tau_B-1} \PE_{X_k}[\tau_B] \right] -  \PE_x[\tau_B] \eqsp.
$$
Therefore, the Markov chain is ergodic of degree 2  if and only if for any $B \in \Xsigma$, $\int_\Xset \pi(dx) \PE_x[\tau_B] < \infty$. The Jensen inequality
(see for instance  \cite[Lemma 3.5]{jarner:roberts:2001}) shows that there exists two positive constants $c_0$ and $b_0$ such that
$P \sqrt{V} \leq \sqrt{V} - c_0 + b_0 \1_C$, and by \cite[Theorem 14.2.3]{meyn:tweedie:1993},  for any $x \in \Xset$, and any $B$ such that $\pi(B) > 0$,
there exists a constant $c(B)$ such that, for ny $x \in \Xset$,
$$
\PE_x[\tau_B] \leq \sqrt{V(x)} + c(B) \eqsp.
$$
Applying to the inequality $PV + c \sqrt{V} \leq V + b \1_C$ shows that $\pi(\sqrt{V}) < \infty$, which concludes the proof.
\end{proof}

We will only give here the scheme of the proof generalizing the approach of \cite {djellout:guillin:2001}, based on the renewal method introduced
(for discrete Markov chains) by Doeblin.
Let us first recall the following crucial result due to Arcones-Ledoux: suppose that $\{ U_i \}$ are \iid\ random variables,
then $b_n^{-1}\sum_{k=1}^n U_k$ satisfies a MDP if and only if
$$\lim_{n\to\infty}{n\over b_n^2}\log\left(n~\PP\left(\|U_k\|\ge b_n\right)\right)=-\infty \eqsp,$$
and the rate function is the natural quadratic one. Note that by an easy approximation argument (at least in the finite dimensional case)
and thus generalizing result by \cite{chen:1997}, the previous condition gives also the MDP for a 1-dependent sequence $\{U_i\}$.

The renewal approach  consists in splitting the sum $S_n \eqdef \sum_{i=0}^{n-1} f(X_i) $ into four different terms:
\begin{equation}
S_n =  \sum_{k=1}^{e(n)}\xi_k + S_{\csigma_0\wedge n}+\left(\sum_{k=1}^{i(n)-1}\xi_k-\sum_{k=1}^{e(n)}\xi_k \right)+\sum_{j=(l(n)+1)}^{n-1}f(X_j)
\label{dec1}
\end{equation}
where $\csigma_0\eqdef\csigma$ and $\csigma_k=\inf\{n>\csigma_{k-1};d_n=1\}$ are the successive return times to the atom of the split chain,
$i(n)\eqdef \sum_{k=0}^{n-1} \1(d_k=1)$ is the number of visits the atom before $n$, $e(n)=\lfloor\epsilon\pi(C)n\rfloor$ is the
expected number of visits to the atom before $n$, $l(n)\eqdef\csigma_{(i(n)-1)\wedge  0} $ is the index of the last visit to the chain to the atom and
$\xi_k\eqdef\sum_{j=\csigma_{k-1}+1}^{\csigma_k}f(X_j)$ is the $f$-modulated moment of the excursion between two successive visits to the atom.

The general idea is to show that only the first term contributes to the moderate deviation principle. To this end we make the following remark:
it can be easily checked that $\{\xi_k\}$ is a sequence of \iid\ random variables
with common distribution
$$
\cPP(\xi_1 \in \cdot)=\cPP_{\check{\nu}}\left(\sum_{j=0}^{\csigma_0}f(X_j) \in \cdot\right) \eqsp.
$$
Note that, when $m > 1$ in (A\ref{assum:smallset}) then the sequence becomes $1$-dependent but essentially the same argument can be carried out.
Under \eqref{eq:theo:mdpunbounded1}-\eqref{eq:theo:mdpunbounded2}, it is easily seen that
\[
\lim_{n\to\infty}{n\over b_n^2}\log\left\{n~\cPP_{\check{\nu}}\left( \left\|\sum_{k=0}^{\check{\sigma}}f(X_k) \right\|\ge b_n\right)\right\}=-\infty \eqsp,
\]
so that the first term satifies a MDP, the identification of the rate function being easily handled.

Consider now the three remaining terms of the \rhs\ of \eqref{dec1}.
We have to show that, for any positive $\varepsilon$
\begin{align}
&\limsup_{n\to\infty}{n\over b_n^2}\log
\cPP_{\check{\mu}}\left(\|S_{\csigma_0\wedge n}\|\ge \varepsilon b_n\right)=-\infty, \label{eq:neg1} \\
&\limsup_{n\to\infty}{n\over b_n^2}\log
\cPP_{\check{\mu}}\left(\left\|\sum_{j=l(n)+1}^{n-1}f(X_j)\right\|\ge \varepsilon b_n\right)=-\infty, \label{eq:neg2} \\
&\limsup_{n\to\infty}{n\over b_n^2}\log \cPP_{\check{\mu}}\left(\left\|\sum_{j=1}^{i(n)-1}\xi_j-\sum_{k=1}^{e(n)}\xi_k\right\|\ge
\varepsilon b_n\right)=-\infty. \label{eq:neg3}
\end{align}

Remark that the condition ensuring the MDP gives directly the first two needed limits. The last one is more delicate, but as seen from the proof done in the atomic case it merely resumes to the MDP of $\left(\check{\sigma}_k-\check{\sigma}_{k-1}-(\epsilon\pi(C))^{-1}\right)$ (given by Arcones-Ledoux result and \eqref{eq:theo:mdpunbounded1}) which enables us to prove that in the sense of moderate deviations the difference $|i(n)-e(n)|$ can be arbitrarily considered of size $\lfloor \delta n\rfloor$ ($\delta$ beeing arbitrary), and the MDP of the sum of $\lfloor \delta n \rfloor$ blocks $(\xi_k)$. This last term being clearly negligible as $\delta$ is arbitrary.
\begin{proof}[Proof of the Theorem \ref{theo:empirMeasurMdp}]
The proof of Theorem \ref{theo:empirMeasurMdp} follows from the projective limit theorem and from the
moderate deviation principle for bounded functions (as stated in Theorem \ref{theo:MDPgen}). The key point consists in checking that
the rate function as expressed in Eq. (\ref{eq:defJf}), Theorem \ref{theo:MDPgen} coincides with the one obtained by the projective limit theorem
(see for instance \cite{acosta:1997} and \cite{acosta:chen:1998}).
\end{proof}

\section{Proof of Theorem \ref{theo:dev}}
\label{sec:mdpdev}
We will the same decomposition than in the
moderate deviations proof, i.e. decomposition (\ref{dec1})
\begin{equation}\label{decdev}
S_n = S_{(\csigma_0)\wedge
n}+\sum_{k=1}^{i(n)-1}\xi_k+\sum_{j=(l(n)+1)}^{n-1}f(X_j).
\end{equation}
We bound $\PP_\mu\left(\left\|\sum_{k=0}^{n-1}f(X_k)\right\|>\varepsilon n\right)$ by $\sum_{i=1}^4 I_i$, where
\begin{align*}
I_1 &\eqdef \cPP_{\check{\mu}}\left(\left\|\sum_{k=0}^{n-1}f(X_k)\right\|>\varepsilon n, \csigma_0>n\right)\\
I_2 &\eqdef \cPP_{\check{\mu}}\left(\left\|S_{\csigma_0} \right\|>{\varepsilon n\over 3}\right)\\
I_3 &\eqdef \cPP_{\check{\mu}}\left(\left\|\sum_{k=1}^{i(n)-1}\xi_k \right\|>{\varepsilon n\over3}\right) \\
I_4 &\eqdef \cPP_{\check{\mu}}\left(\left\| \sum_{j=(l(n)+1)}^{n-1}f(X_j)\right\|>{\varepsilon n\over 3}\right)
\end{align*}
where
$$
I_1 \le  \cPP_{\check{\mu}}\left(\csigma_0>n\right) \le {L\over \Phi^{-1}(n)} \eqsp,
$$
by Theorem \ref{theo:mainresult} if $\mu(V)<\infty$. Remark also
$$
I_2 \le \cPP_{\check{\mu}}\left(\csigma_0>{\varepsilon n\over \|f\|_\infty}\right)
\le {L\over \Phi^{-1}\left({\varepsilon n\over \|f\|_\infty}\right)}\eqsp,
$$
and if $\nu(V)$ is bounded,
\begin{multline*}
I_4 \le \cPP_{\check{\mu}}\left(\max_{k\le /+1}(\csigma_k-\csigma_{k-1})>{\varepsilon n\over \|f\|_\infty}\right)\\
\le (n+1)\cPP_{\check{\nu}}\left(\csigma_0>{\varepsilon n\over \|f\|_\infty}-1\right) \le  {L(n+1)\over \Phi^{-1}\left({\varepsilon n\over \|f\|_\infty}\right)} \eqsp.
\end{multline*}
For the last term, note
\begin{eqnarray*}
I_3 &\le& \cPP_{\check{\mu}}\left(\max_{i\le n}\left\|\sum_{k=1}^{i}\xi_k
\right\|>{\varepsilon n\over 3}\right)\\ 
&\le& 2\cPP_{\check{\mu}}\left(\left\|\sum_{k=1}^{n}\xi_{k}
\right\|>{\varepsilon n\over
6}\right)
\end{eqnarray*}
where the last step follows by Ottaviani's inequality for
i.i.d.r.v. if for $n$ large enough $$\max_{i\le n}
\cPP_{\check{\mu}}\left(\left\|\sum_{k=i}^{n}\xi_{k} \right\|>{\varepsilon n\over
6}\right)\le 1/2.$$

By Chebyschev's inequality, independence and zero mean of the $(\xi_{k})$, it is sufficient to choose $n$ such that

$$n\ge {72\|f\|_\infty^2 \cPE_{\check\nu}((\csigma+1)^2)\over \varepsilon^2},$$
where  $\cPE_{\check\nu}((\csigma+1)^2)$ is finite (and can be easily evaluated) under our drift condition.

By using the Fuk-Nagaev inequality for the
remaining term, we get that for all $y>0$
\begin{eqnarray*}
\cPP_{\check{\mu}}\left(\left\|\sum_{k=1}^{[n/2]+1}\xi_{2k} \right\|>{\varepsilon n\over
12}\right)&\le&([n/2]+1)\cPP_{\check{\mu}}\left(\|\xi_{1}\|>y\right)+\exp\left(-{([n/2]+1)\varepsilon^2\over
(9\cPE\xi_1^2+\varepsilon y)}\right)\\ &\le&{L([n/2]+1)\over
\Phi^{-1}\left(\frac{y}{ \|f\|_\infty}\right)}+\exp\left(-{([n/2]+1)\varepsilon^2\over
(9\cPE\xi_1^2+\epsilon y)}\right),
\end{eqnarray*}
where $\cPE\xi_1^2$ is easily controlled under the drift condition. This concludes the proof.


%

\bibliographystyle{ims}

\end{document}